\documentclass[12pt]{article}
\usepackage{latexsym, amssymb}
\usepackage{mathbbol}

\DeclareMathAlphabet\gothic{U}{euf}{m}{n}

\makeatletter
\def\eqnarray{\stepcounter{equation}\let\@currentlabel=\theequation
\global\@eqnswtrue
\tabskip\@centering\let\\=\@eqncr
$$\halign to \displaywidth\bgroup\hfil\global\@eqcnt\z@
  $\displaystyle\tabskip\z@{##}$&\global\@eqcnt\@ne
  \hfil$\displaystyle{{}##{}}$\hfil
  &\global\@eqcnt\tw@ $\displaystyle{##}$\hfil
  \tabskip\@centering&\llap{##}\tabskip\z@\cr}

\def\endeqnarray{\@@eqncr\egroup
      \global\advance\c@equation\m@ne$$\global\@ignoretrue}
\makeatother

\begin{document}
 \bibliographystyle{tom}

\newtheorem{lemma}{Lemma}[section]
\newtheorem{thm}[lemma]{Theorem}
\newtheorem{cor}[lemma]{Corollary}
\newtheorem{voorb}[lemma]{Example}
\newtheorem{rem}[lemma]{Remark}
\newtheorem{prop}[lemma]{Proposition}
\newtheorem{stat}[lemma]{{\hspace{-5pt}}}

\newenvironment{remarkn}{\begin{rem} \rm}{\end{rem}}
\newenvironment{exam}{\begin{voorb} \rm}{\end{voorb}}

\newcommand{\gota}{\gothic{a}}
\newcommand{\gotb}{\gothic{b}}
\newcommand{\gotc}{\gothic{c}}
\newcommand{\gote}{\gothic{e}}
\newcommand{\gotf}{\gothic{f}}
\newcommand{\gotg}{\gothic{g}}
\newcommand{\gothh}{\gothic{h}}
\newcommand{\gotk}{\gothic{k}}
\newcommand{\gotm}{\gothic{m}}
\newcommand{\gotn}{\gothic{n}}
\newcommand{\gotp}{\gothic{p}}
\newcommand{\gotq}{\gothic{q}}
\newcommand{\gotr}{\gothic{r}}
\newcommand{\gots}{\gothic{s}}
\newcommand{\gotu}{\gothic{u}}
\newcommand{\gotv}{\gothic{v}}
\newcommand{\gotw}{\gothic{w}}
\newcommand{\gotz}{\gothic{z}}
\newcommand{\gotA}{\gothic{A}}
\newcommand{\gotB}{\gothic{B}}
\newcommand{\gotG}{\gothic{G}}
\newcommand{\gotL}{\gothic{L}}
\newcommand{\gotM}{\gothic{M}}
\newcommand{\gotS}{\gothic{S}}
\newcommand{\gotT}{\gothic{T}}

\newcounter{teller}
\renewcommand{\theteller}{\Roman{teller}}
\newenvironment{tabel}{\begin{list}%
{\rm \bf \Roman{teller}.\hfill}{\usecounter{teller} \leftmargin=1.1cm
\labelwidth=1.1cm \labelsep=0cm \parsep=0cm}
                      }{\end{list}}

\newcounter{tellerr}
\renewcommand{\thetellerr}{\roman{tellerr}}
\newenvironment{subtabel}{\begin{list}%
{\rm  \roman{tellerr}.\hfill}{\usecounter{tellerr} \leftmargin=1.1cm
\labelwidth=1.1cm \labelsep=0cm \parsep=0cm}
                         }{\end{list}}

\newcounter{proofstep}
\newcommand{\nextstep}{\refstepcounter{proofstep}\ruimte \par 
          \noindent{\bf Step \theproofstep} \hspace{5pt}}
\newcommand{\firststep}{\setcounter{proofstep}{0}\nextstep}

\newcommand{\Ni}{{\bf N}}
\newcommand{\Ri}{{\bf R}}
\newcommand{\Ci}{{\bf C}}
\newcommand{\Ti}{{\bf T}}
\newcommand{\Zi}{{\bf Z}}
\newcommand{\Fi}{{\bf F}}

\newcommand{\proof}{\mbox{\bf Proof} \hspace{5pt}} 
\newcommand{\remark}{\mbox{\bf Remark} \hspace{5pt}}
\newcommand{\ruimte}{\vskip10.0pt plus 4.0pt minus 6.0pt}

\newcommand{\ad}{{\mathop{\rm ad}}}
\newcommand{\Ad}{{\mathop{\rm Ad}}}
\newcommand{\Aut}{\mathop{\rm Aut}}
\newcommand{\arccot}{\mathop{\rm arccot}}
\newcommand{\diam}{\mathop{\rm diam}}
\newcommand{\divv}{\mathop{\rm div}}
\newcommand{\codim}{\mathop{\rm codim}}
\newcommand{\RRe}{\mathop{\rm Re}}
\newcommand{\IIm}{\mathop{\rm Im}}
\newcommand{\Tr}{{\mathop{\rm Tr}}}
\newcommand{\supp}{\mathop{\rm supp}}
\newcommand{\sgn}{\mathop{\rm sgn}}
\newcommand{\essinf}{\mathop{\rm ess\,inf}}
\newcommand{\esssup}{\mathop{\rm ess\,sup}}
\newcommand{\Int}{\mathop{\rm Int}}
\newcommand{\Leibniz}{\mathop{\rm Leibniz}}
\newcommand{\lcm}{\mathop{\rm lcm}}
\newcommand{\mod}{\mathop{\rm mod}}
\newcommand{\spann}{\mathop{\rm span}}
\newcommand{\ubar}{\underline{\;}}
\newcommand{\one}{\mathbb{1}}

\hyphenation{groups}
\hyphenation{unitary}

\newcommand{\cc}{{\cal C}}
\newcommand{\cd}{{\cal D}}
\newcommand{\ce}{{\cal E}}
\newcommand{\cf}{{\cal F}}
\newcommand{\ch}{{\cal H}}
\newcommand{\ci}{{\cal I}}
\newcommand{\ck}{{\cal K}}
\newcommand{\cl}{{\cal L}}
\newcommand{\cm}{{\cal M}}
\newcommand{\co}{{\cal O}}
\newcommand{\cs}{{\cal S}}
\newcommand{\ct}{{\cal T}}
\newcommand{\cx}{{\cal X}}
\newcommand{\cy}{{\cal Y}}
\newcommand{\cz}{{\cal Z}}

\newfont{\fontcmrten}{cmr10}
\newcommand{\slbrl}{\mbox{\fontcmrten (}}
\newcommand{\slbrr}{\mbox{\fontcmrten )}}

\thispagestyle{empty}

\vspace*{1.5cm}
\begin{center}
{\Large{\bf  Degenerate elliptic operators:}}\\[2mm] 
{\Large{\bf capacity, flux and separation}}  \\[2mm]
\large  Derek W. Robinson$^1$ and Adam Sikora$^2$
\end{center}

\vspace{5mm}

\begin{center}
{\bf Abstract}
\end{center}

\begin{list}{}{\leftmargin=1.8cm \rightmargin=1.8cm \listparindent=10mm 
   \parsep=0pt}
\item

Let $S=\{S_t\}_{t\geq0}$ be the semigroup generated on $L_2(\Ri^d)$ by a self-adjoint, second-order, 
divergence-form, elliptic operator $H$  with Lipschitz continuous coefficients.
Further let $\Omega$ be an open subset of $\Ri^d$ with Lipschitz continuous boundary
$\partial\Omega$.
We prove that $S$ leaves $L_2(\Omega)$ invariant if, and only if, the capacity of 
the boundary with respect to $H$ is zero or if, and only if, the energy flux across the
boundary is zero.
The global result is based on an analogous local result.

\end{list}

\vspace{3cm}
\noindent
July 2005

\vspace{1cm}

\vspace{1cm}

\noindent
{\bf Permanent adresses:}  
\begin{enumerate}
\item
Centre for Mathematics and its Applications\\
Mathematical Sciences Institute\\
Australian National University\\
Canberra, ACT 0200\\
Australia.
\item
 Department of Mathematical Sciences \\ 
  New Mexico State University\\
  P.O. Box 30001\\  
   Las Cruces\\ 
   NM 88003-8001, USA.

\end{enumerate}

\newpage
\setcounter{page}{1}

\section{Introduction}\label{Scd1}

In two earlier papers with Tom ter Elst and Yueping Zhu \cite{ERSZ1} \cite{ERSZ2}
we analyzed the non-ergodic behaviour of degenerate second-order elliptic operators
in divergence form on $\Ri^d$.
In this note we continue the analysis for a class of operators whose coefficients are
Lipschitz continuous.
In particular we establish that the phenomenon of separation highlighted in the earlier
papers can be characterized by the property of zero energy flux across the boundary of 
separation or by zero capacity of the boundary.

Let $h$ be a real, symmetric, bilinear form
\begin{equation}
\varphi,\psi\in C_c^\infty(\Ri^d)\mapsto 
h(\varphi,\psi)=\sum^d_{i,j=1}(\partial_i\varphi,c_{ij}\,\partial_j\psi)
\label{econe1}
\end{equation}
where the coefficients $c_{ij}$ are real $L_\infty$-functions and  the matrix
$C=(c_{ij})$ is symmetric and positive-definite almost-everywhere.
Then the corresponding quadratic form $h(\varphi)=h(\varphi,\varphi)$ is positive.
If the form $h$ is closable  the closure $\overline h$ determines a positive self-adjoint
operator $H$ on $L_2(\Ri^d)$ such that $D(\overline h)=D(H^{1/2})$ and 
$\overline h(\varphi)=\|H^{1/2}\varphi\|_2^2$ (see, for example \cite{Kat1}, Chapter~VI).
The operator $H$ is  interpreted as the second-order elliptic operator with coefficients
$c_{ij}$. 
Since $h(\varphi)\leq \|C\|\,l(\varphi)$ with 
$l(\varphi)=\sum^d_{i=1}\|\partial_i\varphi\|_2^2$ the form of the usual Laplacian $\Delta$
and $\|C\|$ the essential supremum of the norm of the matrix $C(x)$ it follows that 
$W^{1,2}(\Ri^d)\subseteq D(\overline h)$.
Note that this definition does not require any smoothness of the coefficients but does require
closability of the form.

There are two standard settings for the study of second-order elliptic operators.
The Nash--De Giorgi--Aronson theory \cite{Nash} \cite{DG1} \cite{Aro} is based on the strong
ellipticity assumption $C\geq \mu I>0$.
This implies that $h(\varphi)\geq\mu\,l(\varphi)$.
Hence $h$ is closable, $D(\overline h)=W^{1,2}(\Ri^d)$ and the corresponding self-adjoint
operator $H$ satisfies $H\geq\mu\,\Delta$.
The assumption $C\geq \mu I>0$ is a strong condition of non-degeneracy which is equivalent
to the bound $H\geq\mu\,\Delta$ (see, for example, \cite{ERZ1}, Proposition~2).
The principal results of the theory are precise estimates on the global behaviour
of the solutions of the corresponding elliptic and parabolic equations.
In particular the parabolic solutions are bounded above and below by Gaussian functions
expressed in terms of the Riemannian geometry corresponding to the form $h$.

The Fefferman--Phong theory \cite{FP}  is based on the assumption that 
the coefficients $c_{ij}$ are smooth and that $H$ satisfies
the subellipticity  condition
$H\geq \mu\,\Delta^\delta-\nu\,I$ for some $\mu>0$, $\nu\in\Ri$ 
and $\delta\in\langle0,1\rangle$.
The smoothness assumption ensures that $h$ is closable by  Friederichs'  arguments 
(see, for example, \cite{Kat1}, Section~VI.2.3).
Although one still has $W^{1,2}(\Ri^d)\subseteq D(\overline h)$ the domain of the closure
is not readily identifiable.
The advantage of the Fefferman--Phong approach is that it covers several classes of 
degenerate operators.
Moreover, the subellipticity assumption allows the derivation of local versions of many of the estimates
of the strongly elliptic theory
\cite{FSC} \cite{JSC} \cite{JSC1}.
In particular the solutions of the subelliptic parabolic equations still satisfy 
Gaussian upper and lower bounds in terms of the Riemannian geometry.
The major differences between the strongly elliptic and the subelliptic theories are incorporated in the Riemannian geometry.
In the strongly elliptic case the Riemannian distance is equivalent to the Euclidean distance
but this is no longer the case, at least locally, in the subelliptic theory.

In  this note we consider a  somewhat different  situation. 
We mostly assume that the coefficients $c_{ij}$ of $h$ are Lipschitz continuous
but we make no additional ellipticity assumption.
The continuity of the coefficients is sufficient to ensure the form $h$ is closable,
again by Friederichs' arguments.
Moreover, the closure $\overline h$ is a Dirichlet form 
and the corresponding elliptic operator $H$ on $L_2(\Ri^d)$ generates 
a submarkovian semigroup $S$ on the $L_p$-spaces $L_p(\Ri^d)$
(see \cite{FOT} and \cite{BH} for background on Dirichlet forms and submarkovian semigroups).
We also derive some results which are independent of continuity of the coefficients.
In the latter case, however,  there is a difficulty in defining the elliptic operator since the form $h$ is
not necessarily closable (see, for example \cite{FOT}, Section~3.1).
Nevertheless one can consider the strongly elliptic operators
$H_\varepsilon$ associated with the closures of the strongly elliptic forms
$h+\varepsilon \,l$ with $\varepsilon>0$ and then define 
a viscosity operator $H_0$  as the strong 
resolvent limit of the $H_\varepsilon$ as $\varepsilon\to0$.
The viscosity form $h_0$ is then defined by $D(h_0)=D(H_0^{1/2})$ and 
$h_0(\varphi)=\|H_0^{1/2}\varphi\|_2^2$ for $\varphi\in D(h_0)$.
Again it is a Dirichlet form.
A detailed description of the viscosity operator with further background references
can be found in \cite{ERSZ1} and \cite{ERSZ2}.

Our aim is to characterize conditions on the degeneracy of the elliptic operators
which lead to a phenomenon of separation.
In particular we examine conditions for the existence of open subsets $\Omega$ of $\Ri^d$
such that $S_tL_2(\Omega)\subseteq L_2(\Omega)$ for all $t>0$, or in the case of the viscosity 
semigroup $S^{(0)}_tL_2(\Omega)\subseteq L_2(\Omega)$ for all $t>0$.
In Section~\ref{Sldeg2} we prove that the occurrence of separation 
is equivalent to non-ergodicity of  the action of the semigroup
$S$ on $L_2$.
This result does not depend on the detailed structure of the elliptic operator
but is a general  semigroup result.
In the subsequent sections we then derive conditions for separation which are specific
to elliptic operators.
These conditions involve the boundary $\partial\Omega=\overline\Omega\backslash\Omega$
of the set~$\Omega$.
We emphasize that  separation can occur even if the Riemannian distance between 
$\Omega$ and its complement  $\Omega^{\rm c}$ is  finite.

Let $A$ be a general subset of $\Ri^d$ and $k$ positive closed quadratic form 
on $L_2(\Ri^d)$.
The capacity $C_k(A)$ of $A$ with respect to $k$ is defined by
\[
C_k(A)=\inf\{\,C_k(V):V \;{\rm open\; and  }\; V\supset A\,\}
\]
where the capacity $C_k(V)$ of an  open subset $V\subset\Ri^d$ is given by
\[
C_k(V)=\inf\{\, k(\varphi)+\|\varphi\|_2^2:\varphi\in D(k)\,,\, 
\varphi\geq1\;{\rm on }\; V\,\}
\;\;\;.
\]
The general properties of the capacity associated with a Dirichlet form are 
described in \cite{FOT}, Section~2.1, or in \cite{BH}, Section~1.8.

The first result does not require any continuity of the coefficients of the elliptic
operator or of the boundary.

\begin{thm}\label{tcdeg1.0}
Let $h_0$ be the form of the viscosity operator $H_0$ with $L_\infty$-coefficients $c_{ij}$.
Further let $\Omega$ be an open subset of $\Ri^d$.

If $C_{h_0}(\partial\Omega)=0$ then $S^{(0)}_tL_2(\Omega)\subseteq L_2(\Omega)$.
\end{thm}

The proof of this result will be given in Section~\ref{Scdeg2.0} together with an
example which shows that a strict converse is not correct.
One can have separation without the capacity of the boundary being zero.
Nevertheless the next theorem establishes that if the coefficients and the boundary 
are both Lipschitz continuous then the converse to Theorem~\ref{tcdeg1.0} is valid.
The theorem also gives a characterization of separation in terms of energy flux.

If $A$ is a Lipschitz continuous hypersurface then there is a unique (up to orientation)
normal vector  $n_A(x)=(n_{A,1}(x),\ldots,n_{A,d}(x))$ at almost every $x\in A$.
If the orientations are  chosen in a  consistent manner then the (energy) flux across
$A$ corresponding to $ h$ is defined by $F_h(A)=(F_h(A)_1,\ldots, F_h(A)_d)$
where $F_h(A)_i=\sum^d_{j=1}c_{ij}\,n_{A,j}$.
Since $C=(c_{ij})\geq0$ it  follows that  $F_h(A)=0$ on $A$ if, and only if,
\[
(n_A,Cn_A)=\sum^d_{i,j=1}c_{ij}\,n_{A,i}\,n_{A,j}=0 
\]
on $A$. 
But this is in turn equivalent to the normal vector $n_A$ being an eigenvector 
of the matrix $C$ corresponding to the eigenvalue zero.

\begin{thm}\label{tcdeg1.2}
Assume that the coefficients $c_{ij}$ of $h$ 
and  the boundary $\partial\Omega $ of the  open subset
 $\Omega\subseteq \Ri^d$ are both Lipschitz continuous.

\smallskip

The following conditions are equivalent.
\begin{tabel}
\item
\label{tcdeg1.2-1}
$\;S_tL_2(\Omega)\subseteq L_2(\Omega)$ for all $t\geq 0$.
\item
\label{tcdeg1.2-2}
$\;C_h(\partial\Omega)=0$.
\item
\label{tcdeg1.2-3}
$\;F_h(\partial\Omega)=0$ almost everywhere on $\partial\Omega$.
\end{tabel}
\end{thm}

In fact the equivalence of zero capacity of the boundary and zero flux across 
the boundary is valid in a broader context.
It is not necessary that the hypersurface separates $\Ri^d$ into disjoint components.

\begin{thm}\label{tcdeg1.1}
Assume that the coefficients $c_{ij}$ of $h$ 
are  Lipschitz continuous.
Further let $A$ be a Lipschitz continuous hypersurface in~$\Ri^d$.

\smallskip

The following conditions are equivalent.
\begin{tabel}
\item
\label{tcdeg1.1-1}
$\;C_h(A)=0$.
\item
\label{tcdeg1.1-2}
$\;F_h(A)=0$ almost everywhere on $A$.
\end{tabel}
\end{thm}

The proof of Theorem~\ref{tcdeg1.1} is based on a property of local separation.
It is given in Section~\ref{Scdeg2} together with a proof of 
 Theorem~\ref{tcdeg1.2}.
Then in Section~\ref{Scdeg3} we briefly discuss various possible extensions
and related properties.

\section{Ergodicity and separation}\label{Sldeg2}

In this section we establish general criteria for separation in terms
of the action of the semigroup.
Proposition~\ref{perg} is  similar to well known extensions
of the Perron--Frobenius theorem (see, for example, \cite{RS4}, Section~XIII.12)
but we do not require the existence of any  point spectrum.
The main ingredient is positivity of the semigroup.
The result is derived for a positive, self-adjoint, strongly continuous
semigroup  $S=\{S_t\}_{t\geq 0}$ acting on the Hilbert space $L_2(X)\,(=L_2(X\,;\mu))$
where  $(X,\mu)$ is  a $\sigma$-finite measure space.

The semigroup $S$ is defined to be {\bf ergodic} if $(\varphi, S_t\psi)>0$ for 
each pair of non-zero, non-negative $\varphi,\psi\in L_2(X)$ and at least one $t>0$.
Similarly $S$ is defined to be {\bf strictly positive} if $(\varphi, S_t\psi)>0$ for 
each pair of non-zero, non-negative $\varphi,\psi\in L_2(X)$ and all $t>0$.
Further a family of operators acting on $L_2(X)$ is defined to be {\bf irreducible}
if there is no non-trivial closed subspace of $L_2(X)$ which is left invariant by the action of 
the family.
Note that the  bounded measurable functions $L_\infty(X)$ act as multipliers on $L_2(X)$.
The phenomenon of separation is closely tied  to the reducibility of the family $S\cup L_\infty$
formed by the operators $S_t$, $t>0$, together with the multiplication 
operators $L_\infty(X)$.
\begin{prop}\label{perg}
The following conditions are equivalent.
\begin{tabel}
\item\label{perg1}
$S$ is ergodic.
\item\label{perg2}
$S$ is strictly positive.
\item\label{perg3}
$S\cup L_\infty$ is irreducible.\end{tabel}
\end{prop}
\proof\
It is evident that \ref{perg2}$\Rightarrow$\ref{perg1} but 
\ref{perg1}$\Rightarrow$\ref{perg2} by \cite{RS4}, Theorem~XIII.44 (see also 
\cite{ERSZ2}, Lemma~4.1).

\smallskip

\noindent\ref{perg1}$\Rightarrow$\ref{perg3}$\;$
Assume \ref{perg3} is false. 
Let $\ch$ denote a non-trivial closed subspace which is invariant under the action
of the operators $S\cup L_\infty$.
If $\psi\in\ch$ then $|\psi|\in\ch$ because $|\psi|=\sgn\psi\cdot\psi$ 
and $\sgn\psi\in L_\infty$.
The orthogonal complement $\ch^\perp$ is alsoa non-trivial closed subspace which is
invariant under the $S\cup L_\infty$.
Hence if $\varphi\in \ch^\perp$ then $|\varphi|\in\ch^\perp$.
But then $(|\varphi|, S_t |\psi|)=0$ for
all $t>0$ because $\ch$ is $S$-invariant.
Hence \ref{perg1} is false.

\smallskip

\noindent\ref{perg3}$\Rightarrow$\ref{perg1}$\;$
Assume \ref{perg1} is false.
Then there are 
non-negative, non-zero, $\varphi,\psi\in L_2(X)$ such that $(\varphi,S_t\psi)=0$
for all $t>0$.
Let $\ck$ be the convex cone in   $L_2(X)$ spanned by the non-negative $\chi$
such that $(\chi,S_t\psi)=0$ for all $t>0$.
Now $\varphi\in \ck$ so $\ck$ is non-empty.
But $(\psi,S_t\psi)=\|S_{t/2}\psi\|_2^2>0$ since $\psi\neq0$.
Thus $\psi\not\in\ck$.
Clearly $\ck$ is closed and invariant under the action of the semigroup $S$.
Next we argue that it is also invariant under multiplication by non-negative 
$\eta\in L_\infty(X)$.

Each  non-negative $\eta\in L_\infty(X)$
 can be approximated monotonically from below by simple functions
$\sum^n_{i=1}\eta_i\one_{A_i}$ where the $\eta_i>0$ and 
the $A_i$ are measurable subsets.
Therefore to prove that $\eta\,\ck\subseteq \ck$ it  suffices to prove that 
$(\one_A\chi,S_t\psi)=0$ for all  $\chi\in \ck$, all $t>0$ and each measurable subset $A$.
Again $\chi$ can be approximated monotonically from below by simple functions
$\sum^n_{i=1}\chi_i\one_{A_i}$ with  $\chi_i>0$.
Then, since $S$ is positive and $\psi$ is non-negative,
\[
\chi_i\,(\one_{A_i},S_t\psi)\leq (\sum^n_{i=1}\chi_i\one_{A_i},S_t\psi)
\leq (\chi,S_t\psi)=0
\]
for all $t>0$.
Therefore $(\one_{A_i},S_t\psi)=0$ for all $t>0$.
Then, however, 
\[
(\one_A\one_{A_i},S_t\psi)=(\one_{A\cap A_i}, S_t\psi)\leq (\one_{A_i}, S_t\psi)=0
\]
for all $t>0$.
Hence
\[
(\one_A\sum^n_{i=1}\chi_i\one_{A_i},S_t\psi)=0
\]
for all $t>0$.
Taking a limit of the approximants one then concludes that  $(\one_A\,\chi, S_t\psi)=0$
for all $t>0$.
Hence  $\one_A\,\chi\in\ck$.

Next introduce the  subspace $\ch=\ck-\ck$ of $L_2(X)$.
It follows from the invariance properties of $\ck$ that $\ch$ is invariant under
the action of $S\cup L_\infty$.
In particular if $\chi\in\ch$ then $|\chi|=\sgn\chi\cdot\chi\in\ch$.
But it then follows by the definition of $\ck$ that $(\chi, S_t\psi)=0=(|\chi|, S_t\psi)$
for all $t>0$.
Therefore $\chi_\pm=|\chi|\pm\chi\in \ck$.
Now suppose that the sequence $\chi_n\in \ch$ converges to $\chi\in L_2$.
Then $|\chi_n|\in\ch$ converges to $|\chi|\in L_2$ and $|\chi_n|\pm\chi_n\in\ck$ 
converges to $\chi_\pm=|\chi|\pm\chi\in L_2$.
Since $\ck$ is closed it follows that $\chi_\pm\in\ck$ and $\chi=\chi_+-\chi_-\in\ch$.
This proves that $\ch$ is closed.
But $\varphi\in \ch$.
Moreover, $\psi\geq0$ and $\psi\not\in\ck$.
Hence  $\psi\not\in\ch$.
Therefore $\ch$ is non-trivial.

In summary $\ch$ is a non-trivial closed subspace of $L_2(X)$ which is invariant 
under the action of $S\cup L_\infty$. 
Hence  Condition~\ref{perg3} is false and so \ref{perg3}$\Rightarrow$\ref{perg1}.
\hfill$\Box$

\ruimte

One can immediately characterize the separation property by failure of the equivalent
conditions of Proposition~\ref{perg}.

\begin{cor}\label{cerg} The following conditions are equivalent.
\begin{tabel}
\item\label{cerg1}
There exists a measurable set $A$  such that  $L_2(A)$ is a non-trivial subspace
of $L_2(X)$ and $S_tL_2(A)\subseteq L_2(A)$ for all $t>0$.
\item\label{cerg2}
$S$ is not ergodic.
\end{tabel}
\end{cor}
\proof\
\ref{cerg1}$\Rightarrow$\ref{cerg2}$\;$
If $\varphi\in L_2(A)$ and $\psi\in L_2(A^{\rm c})=L_2(A)^\perp$ then $(\varphi, S_t\psi)=0$
for all $t>0$. 
Hence $S$ is not ergodic.

\smallskip

\noindent\ref{cerg2}$\Rightarrow$\ref{cerg1}$\;$
Since $S$ is not ergodic it follows from Proposition~\ref{perg} that 
$S\cup L_\infty$ is not irreducible.
Let $\ch$ be  a non-trivial closed subspace which is invariant under $S\cup L_\infty$.
Further let $E$ denote  the orthogonal projection onto $\ch$.
Then $E$ commutes with each $S_t$ and with $L_\infty(X)$.
But the algebra of multipliers $L_\infty(X)$ is a maximal abelian von Neumann algebra.
Therefore $L_\infty(X)$ is its own commutant. 
Hence  $E\in L_\infty(X)$.
Now set $A=\{\,x: E(x)>0\,\}$.
Then $\ch=L_2(A)$.
\hfill$\Box$

\ruimte

Note that since the commutant of  the family $S\cup L_\infty$ is a closed $^*$subalgebra
of $L_\infty$ it is abelian.
Hence $S\cup L_\infty$ can be uniquely decomposed into irreducible components onto
 subspaces $L_2(A)$ of $L_2(X)$.

\section{Capacity and separation}\label{Scdeg2.0}

   In this section we prove Theorem~\ref{tcdeg1.0} and give a counterexample
to its converse. 
One difficulty is that the theorem is formulated in terms of the Dirichlet form $h_0$ 
associated with the viscosity operator $H_0$ and this is linked to the original form 
$h$ in a rather indirect manner.
Therefore one needs to use some of the general structure of local Dirichlet forms in
its proof.
We continue to work in the general context of Section~\ref{Sldeg2} and 
follow the arguments of \cite{ERSZ2}.

Let  $\ce$ be a Dirichlet form on $L_2(X)$.
Then   for all $\varphi \in D(\ce) \cap L_\infty(X)$ define the truncated form
$\ci^{(\ce)}_\varphi \colon D(\ce) \cap L_\infty(X) \to \Ri$ by
\[
\ci^{(\ce)}_\varphi(\psi) = \ce(\varphi \, \psi,\varphi) - 2^{-1} \ce(\psi,\varphi^2)
 \]
where all functions are real-valued.
If $\psi\geq 0$ it follows that 
$\varphi\mapsto \ci^{(\ce)}_\varphi(\psi)\in \Ri$
is a Markovian form with  domain $D(\ce)\cap L_\infty(X)$
(see  \cite{BH}, Proposition~I.4.1.1).
Secondly,  define
\[
|||\ci^{(\ce)}_\varphi|||
= \sup \{ \,| \ci^{(\ce)}_\varphi(\psi)| : 
      \psi \in D(\ce) \cap L_{\infty}(X) , \; \|\psi\|_1 \leq 1\, \}
\in [0,\infty]
\;\;\; .  
\]
The form
$\ce$ is defined to be regular if there is a subset of $D(\ce) \cap C_c(X)$
which is  a core of $\ce$, i.e., which is dense  in $D(\ce)$
with respect to the  norm $\varphi \mapsto (\ce(\varphi) + \|\varphi\|_2^2)^{1/2}$,
and which is also dense  in $C_0(X)$ with respect to the supremum norm 
(see \cite{FOT}, Section~1.1).
Moreover, the form is defined to be local if $\ce(\psi,\varphi) = 0$ for all 
$\varphi,\psi \in D(\ce)$
and $a \in \Ri$ such that $(\varphi + a \one) \psi = 0$
(see \cite{BH}, Section~I.5).

Then the following statement is a simplified version of Lemma~3.4 in \cite{ERSZ2}.

\begin{lemma} \label{lcdeg2.01}
Let $\ce$ be a local, regular, Dirichlet form.
If $\varphi \in D(\ce)$ and $\psi \in D(\ce) \cap L_\infty(X)$
with $||| \ci^{(\ce)}_\psi||| < \infty$ then 
$\psi \, \varphi \in D(\ce)$ and 
\[
\ce(\psi \, \varphi)^{1/2}
\leq ||| \ci^{(\ce)}_\psi|||^{1/2} \, \|\varphi\|_2
   + \|\psi\|_\infty \, \ce(\varphi)^{1/2}
\;\;\; .  \]
\end{lemma}

Now we turn to the specific context of the elliptic form $h$ on $L_2(\Ri^d)$ introduced
in Section~\ref{Scd1}.
The lemma then  applies directly to the  viscosity form $h_0$ since it is regular by Lemma~2.1 
of \cite{ERSZ2} and local by Proposition~2.2 of \cite{ERSZ2}.

The  important feature of the proof of Theorem~\ref{tcdeg1.0} is the observation
 that sets with capacity zero can effectively be neglected.

\begin{prop}\label{pcdeg2.1}
Let $A$ be a subset of $\Ri^d$.
The following conditions are equivalent.

\begin{tabel}
\item\label{pcdeg2.1-1} $\;\;\;$ $C_{h_0}(A)=0$.
\item\label{pcdeg2.1-2}
$\;\;\;$ $D(h_0)\cap L_{2,c}(\Ri^d\backslash A)$ is a core of $h_0$.
\end{tabel}
\end{prop}
\proof\ 
\ref{pcdeg2.1-1}$\Rightarrow$\ref{pcdeg2.1-2}$\;$
First note that $L_{2,c}(\Omega)$ is defined for any subset $\Omega$ of $\Ri^d$
as the functions in $L_2(\Ri^d)$ with compact support in $\Omega$.
Secondly, note that $C_c^\infty(\Ri^d)$ is a core of $h_0$ by Lemma~2.1 of \cite{ERSZ2}.
Therefore it suffices to prove that each $\varphi\in C_c^\infty(\Ri^d)$ can be approximated
by a sequence $\varphi_n\in D(h_0)\cap L_{2,c}(\Ri^d\backslash A)$
with respect to the  norm $\varphi \mapsto (h_0(\varphi) + \|\varphi\|_2^2)^{1/2}$.

  Fix $\varphi\in C_c^\infty(\Ri^d)$.
Since $C_{h_0}(A)=0$ one may choose  sequences $\chi_n\in D(h_0)$ and open subsets 
$U_n\supset A$ such that $h_0(\chi_n)+\|\chi_n\|_2^2\leq n^{-1}$ and $\chi_n\geq 1$ on $U_n$.
But it follows from the Dirichlet property that if $\chi_n\in D(h_0)$ then 
$\chi_n\wedge 1\in D(h_0)$.
Moreover, $h_0(\chi_n\wedge 1)\leq h_0(\chi_n)$ and $\|\chi_n\wedge 1\|_2\leq \|\chi_n\|$.
Therefore replacing $\chi_n$ by $\chi_n\wedge 1$ one may assume
$h_0(\chi_n)+\|\chi_n\|_2^2\leq n^{-1}$ and $\chi_n= 1$ on $U_n$.
Now set $\varphi_n=\varphi(1-\chi_n)$.
Then $\supp\varphi_n\subseteq \Ri^d\backslash U_n$.
But 
\[
\|\varphi-\varphi_n\|_2^2=\|\varphi\chi_n\|_2^2
\leq \|\varphi\|_\infty^2\|\chi_n\|_2^2\leq n^{-1}\|\varphi\|_\infty^2
\;\;\;.
\]
Moreover, by Lemma~\ref{lcdeg2.01} applied to $h_0$ one has
\[
h_0(\varphi-\varphi_n)=h_0(\varphi\chi_n)
\leq 2\,||| \ci^{(h_0)}_\varphi|||\cdot \|\chi_n\|_2^2 +2\,\|\varphi\|_\infty^2\,h_0(\chi_n)
\;\;\;.
\]
Now we must estimate $||| \ci^{(h_0)}_\varphi|||$.

The map $\ce\mapsto  \ci^{(\ce)}_\varphi$ is monotonic by Proposition~3.2 of \cite{ERSZ2}
and $h_0\leq h\leq \|C\|\,l$.
Therefore $||| \ci^{(h_0)}_\varphi|||\leq \|C\|\cdot||| \ci^{(l)}_\varphi|||$.
But one calculates immediately that
\[
\ci^{(l)}_\varphi(\psi)=\sum^d_{i=1}(\partial_i\varphi,\psi\partial_i\varphi)
\]
and so $||| \ci^{(l)}_\varphi|||\leq \|\varphi\|_{1,\infty}^2
=\sup_{x\in\Ri^d}\sum^d_{i=1}(\partial_i\varphi)(x)^2$.
Combining these estimates one has
\[
h_0(\varphi-\varphi_n)\leq 
2\,\|\varphi\|_{1,\infty}^2\,\|\chi_n\|_2^2 +2\,\|\varphi\|_\infty^2\,h_0(\chi_n)
\leq 2\,n^{-1}(\|\varphi\|_{1,\infty}^2+\|\varphi\|_\infty^2)
\;\;\;.
\]
Thus $h_0(\varphi-\varphi_n)+\|\varphi-\varphi_n\|_2^2\to 0$ as $n\to\infty$.

\smallskip

\noindent\ref{pcdeg2.1-2}$\Rightarrow$\ref{pcdeg2.1-1}$\;$
It follows from  monotonicity of the capacity that it suffices to prove $C_{h_0}(A_n)$
for each bounded subset $A_n$ of $A$.
Thus we can effectively assume that $A$ is bounded. 
Then let   $\varphi\in D(h_0)$ and $\varphi\geq 1$ on an open neighbourhood $U$ of $A$.
By hypothesis there exists a sequence $\varphi_n\in D(h_0)\cap L_{2,c}(\Ri^d\backslash A)$ 
such that $h_0(\varphi-\varphi_n)+\|\varphi-\varphi_n\|_2^2\to 0$ as $n\to\infty$.
Since  $\varphi_n$ has compact support in $\Ri^d\backslash A$ it also follows that there
is an  open neighbourhood $U_n$ of $A$ such that $\varphi_n=0$ on $U_n$.
Therefore $\varphi-\varphi_n\geq 1$ on $U\cap U_n$ and one must have $C_{h_0}(A)=0$.
\hfill$\Box$

\ruimte

Now the proof of Theorem~\ref{tcdeg1.0} is straightforward.

\smallskip

\noindent{\bf Proof of Theorem~\ref{tcdeg1.0}}$\;$
The invariance property $S^{(0)}_tL_2(\Omega)\subseteq L_2(\Omega)$ for all $t>0$ 
is equivalent to the statement that for each $\varphi\in D(h_0)$ one has $\varphi\one_\Omega\in D(h_0)$
and
\[
h_0(\varphi)=h_0(\varphi\one_\Omega)+h_0(\varphi\one_{\Omega^{\rm c}})
\]
where $\one_\Omega$ denotes the characteristic function of the set $\Omega$.
The equivalence of these properties is established in Lemma~6.3 of \cite{ERSZ1}.
Therefore it suffices to establish this decomposition property.

If $\varphi\in D(h_0)$ it follows from Proposition~\ref{pcdeg2.1} that there exists
a sequence $\varphi_n\in D(h_0)\cap L_{2,c}(\Ri^d\backslash\partial\Omega)$
such that $h_0(\varphi-\varphi_n)+\|\varphi-\varphi_n\|_2^2\to0$ as $n\to\infty$.
Then each $\varphi_n$ has a unique decomposition 
$\varphi_n=\varphi_n^{(1)}+\varphi_n^{(2)}$ into functions with compact
support $\supp\varphi_n^{(1)}\subseteq \Omega$ and 
$\supp\varphi_n^{(2)}\subseteq \Omega^{\rm c}$.
Moreover, one may choose $\chi\in C_c^\infty(\Ri^d)$ such that 
$\varphi_n^{(1)}=\varphi_n\chi$.
Therefore $\varphi_n^{(1)},\varphi_n^{(2)}\in D(h_0)$.
But $\varphi_n^{(1)}=\varphi\one_\Omega$ and $\varphi_n^{(2)}
=\varphi\one_{\Omega^{\rm c}}$.
Thus $\varphi_n\one_\Omega,\varphi_n\one_{\Omega^{\rm c}}\in D(h_0)$.
Moreover $(\varphi_n\one_\Omega,\varphi_n\one_{\Omega^{\rm c}})=0$.
Hence $h_0(\varphi_n\one_\Omega,\varphi_n\one_{\Omega^{\rm c}})=0$ by locality.
Thus 
\[
h_0(\varphi_n)=h_0(\varphi_n\one_\Omega)+h_0(\varphi_n\one_{\Omega^{\rm c}})
\;\;\;.
\]
But one can make a similar argument with $\varphi_n$ replaced by $\varphi_n-\varphi_m$
to obtain 
\[
h_0(\varphi_n-\varphi_m)=h_0((\varphi_n-\varphi_m)\one_\Omega)+
h_0((\varphi_n-\varphi_m)\one_{\Omega^{\rm c}})
\;\;\;.
\]
The decomposition property for $\varphi$ then follows by continuity.
\hfill$\Box$

\ruimte

Next we demonstrate by example that a strict converse of the theorem is not valid; 
separation does not automatically imply that the boundary has zero capacity.

\begin{exam}\label{excdeg2.1}
Let $d=1$. Then $h(\varphi)=\int_\Ri c\,(\varphi')^2$ with $c\geq0$.
We choose $c$ such that $c(x)=(x^2/(1+x^2))^{1/2}$ if $x\geq0$ and 
$c(x)=(x^2/(1+x^2))^\delta$  if $x<0$ where $\delta\in\langle0,1/2\rangle$.
Then $h$ is closable and the system separates into two subsystems on the right
and left half lines (see \cite{ERSZ1}, Propositions~2.3 and 6.5).
The separation occurs because $c(x)=0(x)$ as $x\to0_+$. 
The coefficient is continuous but it is not Lipschitz continuous because
$\delta<1/2$.
The boundary of separation is the point $\{0\}$ and we argue that $C_{\overline h}(\{0\})>0$.

Let $\varphi\in W^{1,2}(\Ri)$ with $\varphi=1$ on an interval 
$\langle-\varepsilon, \varepsilon\rangle$.
Define $\varphi_+$ by $\varphi_+(x)=\varphi(x)$ if $x<0$ and $\varphi_+(x)=\varphi(-x)$ if
$x\geq0$.
Then $\varphi_+\in W^{1,2}(\Ri)$ and  $\varphi_+=1$ on 
$\langle-\varepsilon, \varepsilon\rangle$.
But
\begin{eqnarray*}
h(\varphi)+\|\varphi\|_2^2&=&
\int_{-\infty}^\infty dx\,\Big(c(x)\,(\varphi'(x))^2+\varphi(x)^2\Big)\\[5pt]
&\geq&\int_{-\infty}^0 dx\,\Big(c(x)\,(\varphi'(x))^2+\varphi(x)^2\Big)
=2^{-1}\int_{-\infty}^\infty dx\,\Big(c_\delta(x)\,(\varphi_+'(x))^2+\varphi_+(x)^2\Big)
\end{eqnarray*}
where $c_\delta(x)=(x^2/(1+x^2))^\delta$ for all $x\in \Ri$.
Thus if $h_\delta$ denotes the form with coefficient $c_\delta$ one has
\[
h(\varphi)+\|\varphi\|_2^2\geq 
2^{-1}\Big(h_\delta(\varphi_+)+\|\varphi_+\|_2^2\Big)
\;\;\;.
\]
The form $h_\delta$ is closable by \cite{ERSZ1}, Proposition~2.3.
Then since $W^{1,2}(\Ri)$ is a core of $h$  it follows that for each 
$\varphi\in D(\overline h)$ with $\varphi=1$ on 
$\langle-\varepsilon, \varepsilon\rangle$ there is a $\varphi_+\in D(\overline h_\delta)$
with $\varphi_+=1$ on 
$\langle-\varepsilon, \varepsilon\rangle$ such that 
\[
\overline h(\varphi)+\|\varphi\|_2^2\geq 
2^{-1}\Big(\overline h_\delta(\varphi_+)+\|\varphi_+\|_2^2\Big)
\;\;\;.
\]
Therefore 
\[
C_{\overline h}(\{0\})\geq 2^{-1}C_{\overline h_\delta}(\{0\})
\;\;\;.
\]
Next let $H_\delta$ denote the self-adjoint operator associated with 
$\overline h_\delta$ and $\Delta$ the self-adjoint version of $-d^2/dx^2$ on $L_2(\Ri)$.
Then it follows from  Example~5.6 in \cite{ERSZ1} that one has the subellipticity estimate
$I+H_\delta\geq \omega\,(I+\Delta)^{1-\delta}$ for some $\omega>0$. 
Therefore
\[
C_{\overline h}(\{0\})\geq 
(\omega/2)\inf\{\,(\varphi,(I+\Delta)^{1-\delta}\varphi):\varphi(\{0\})=1\}
\;\;\;.
\]
The infimum can be easily calculated by Fourier transformation and one concludes that 
\[
C_{\overline h}(\{0\})\geq 
(\omega/2)\Big(\int_\Ri dp\,(1+p^2)^{-(1-\delta)}\Big)^{-1}>0
\;\;\;.
\]
Note that the bound is finite since $\delta<1/2$ but it tends to zero as $\delta\to1/2$.
\end{exam}

The separation properties of the one-dimensional example with the form $h_\delta$
 can be understood in terms of  a related distance.
Theorem~1.3 in \cite{ERSZ2} characterizes separation in terms of a set-theoretic `distance'
which is defined by a variational principle with  trial functions in 
$D(\overline h_\delta)\cap L_\infty$.
The system separates into the  components $\langle-\infty,0\rangle$ and 
$\langle 0,\infty\rangle$ if and only if the distance between 
the left and right is infinite.
But in the example the subellipticity condition 
$I+H_\delta\geq \omega(I+\Delta)^{1-\delta}$ ensures that 
$D(\overline h_\delta)\subseteq D((I+\Delta)^{(1-\delta)/2})\subseteq C_0(\Ri)$
if $\delta<1/2$.
This  continuity property is sufficient to imply that the set-theoretic distance 
coincides with the corresponding Riemannian distance.
Hence the distance cannot be infinite and separation cannot take place. 
This gives an indirect confirmation that $C_{\overline h_\delta}(\{0\})>0$.

It should also be emphasized that one can have separation even if the Riemannian distance
is finite. 
In the one-dimensional example with the form $h_\delta$ one has separation for
$\delta\in[1/2,1\rangle$ but the Riemannian distance between $x$ and $y$ is 
$d(x\,;y)=|\int^y_x c_\delta^{-1/2}|$ which is finite for all $x,y$.
Thus the evolution described by the semigroup cannot penetrate from $x<0$ to $y>0$ although
$d(x\,;y)<\infty$.

The subelliptic estimate for the capacity derived in the one-dimensional example extends
to higher dimensions at least for bounded open subsets.
Assume the subellipticity estimate
\[
h(\varphi)\geq \mu\,\|\Delta^{(1-\delta)/2}\varphi\|_2^2-\nu\,\|\varphi\|_2^2
\]
is valid for some $\mu>0,\nu\geq 0$, $\delta\in[0,1\rangle$  and all $\varphi\in C_c^\infty(\Ri^d)$.
Then a similar estimate is true for the viscosity form $h_0$.
Hence one has an estimate
\[
I+H_0\geq \omega_\delta\,(I+\Delta)^{1-\delta}
\]
in the sense of quadratic forms for some $\omega_\delta>0$.
Therefore if $U\subset \Ri^d$ is an open set and $|U|<\infty$ one has
\begin{eqnarray*}
C_{h_0}(U)&\geq&
\omega_\delta\inf\{\,\|(I+\Delta)^{(1-\delta)/2}\psi\|_2^2:
\,\psi\in D(\Delta^{(1-\delta)/2})\,,\,\psi=1\,{\rm\, on\;} U\}\\[5pt]
&\geq&
\omega_\delta\inf\{\,\|(I+\Delta)^{(1-\delta)/2}\psi\|_2^2:
\,\psi\in D(\Delta^{(1-\delta)/2})\,,\,(\psi,\one_U)=|U|\,\}
\;\;\;.
\end{eqnarray*}
But the last infimum is readily calculated. 
One obtains
\begin{equation}
C_{h_0}(U)\geq \omega_\delta\,|U|^2\,(\one_U,(I+\Delta)^{-(1-\delta)}\one_U)^{-1}
\;\;\;.
\label{ecdeg2.73}
\end{equation}
The infimum is attained with $\psi=|U|\,\Delta^{-(1-\delta)}\one_U/(\one_U,(I+\Delta)^{-(1-\delta)}\one_U)$.
We will use this estimate in the next section to establish that separation with a Lipschitz continuous 
surface is not possible if $\delta<1/2$.

Finally we note that it is unclear whether the condition $C_h(\partial\Omega)=0$ implies separation if the 
coefficients are Lipschitz continuous but the boundary is not.
The example shows that continuity of the boundary without Lipschitz continuity of the 
coefficients does not suffice for the implication.

\section{{Capacity and flux}}\label{Scdeg2}

In this section we give the proofs of  Theorems~\ref{tcdeg1.2} and \ref{tcdeg1.1}.
Therefore we assume throughout the section that the coefficients are Lipschitz continuous
which then implies that the form $h$ is closable.

The implication \ref{tcdeg1.2-2}$\Rightarrow$\ref{tcdeg1.2-1} in Theorem~\ref{tcdeg1.2}
is a corollary of Theorem~\ref{tcdeg1.0} which was established in Section~\ref{Scdeg2.0}.
Next we prove that \ref{tcdeg1.2-1}$\Rightarrow$\ref{tcdeg1.2-3} 
in Theorem~\ref{tcdeg1.2}.
Then it remains to prove that \ref{tcdeg1.2-3}$\Rightarrow$\ref{tcdeg1.2-2}.
But this will be a corollary of   Theorem~\ref{tcdeg1.1}

\smallskip

\noindent{\bf Proof of \ref{tcdeg1.2-1}$\Rightarrow$\ref{tcdeg1.2-3}  in Theorem~\ref{tcdeg1.2}}$\;$ 
First to avoid confusion let $S^{(p)}$ denote the submarkovian semigroup $S$ acting on $L_p$.
then it follows from \cite{ERSZ1},
Lemma~6.1, that the separation property of Condition~\ref{tcdeg1.2-1} is equivalent to the conservation property 
$S^{(\infty)}_t\one_\Omega=\one_\Omega$ for all $t\geq0$.
Now   if $\varphi\in C_c^\infty(\Ri^d)$ then
\[
(\one_\Omega,\varphi)=
(S^{(\infty)}_t\one_\Omega,\varphi)=
(\one_\Omega,S^{(1)}_t\varphi)=
(\one_\Omega,S^{(2)}_t\varphi)
\]
for all $t\geq0$.
Hence
\[
0=(\one_\Omega,\varphi)-(\one_\Omega,S^{(2)}_t\varphi)
=\int^t_0ds\,(\one_\Omega,S^{(2)}_sH\varphi)
=t\,(\one_\Omega,H\varphi)
\]
for all $t>0$.
But then
\[
0=(\one_\Omega,H\varphi)=
\sum^d_{i,j=1}\int_\Omega dx\,(\partial_i\,c_{ij}\partial_j\,\varphi)(x)
=\int_\Omega dx\,\divv\Psi
\]
where $\Psi=\sum^d_{j=1}c_{ij}\,\partial_j\varphi$.
Therefore it follows from   Stokes' theorem that
$\sum^d_{i=1}n_{\Omega,i}\,c_{ij}=0$ on $\partial\Omega$.
\hfill$\Box$

\bigskip

It now remains to prove Theorem~\ref{tcdeg1.1}.

It is convenient to introduce the
  energy density as a symmetric bilinear form over $C_c^\infty(\Ri^d)\times C_c^\infty(\Ri^d)$
with values in $L_1(\Ri^d)$ by
\begin{equation}
\Gamma_{\varphi,\psi}=\sum^d_{i,j=1}c_{ij}\,(\partial_i\varphi)\,(\partial_j\psi)
\;\;\;.
\label{econe2}
\end{equation}
Then 
\[
h(\varphi,\psi)=\int_{\Ri^d}dx\,\Gamma_{\varphi,\psi}(x)
\;\;\;.
\]
The corresponding quadratic form $\varphi\in C_c^\infty(\Ri^d)\mapsto \Gamma_\varphi=\Gamma_{\varphi,\varphi}$
is positive and the truncated form $\ci^{(\overline h)}$ defined in Section~\ref{Scdeg2.0}  is given by
\[
\ci^{(\overline h)}_\varphi(\psi) = \overline h(\varphi \, \psi,\varphi) - 2^{-1} \overline h(\psi,\varphi^2)
=\int_{\Ri^d}dx\,\Gamma_\varphi(x)\,\psi(x)
\]
for all $\varphi,\psi\in  C_c^\infty(\Ri^d)$.
It follows from positivity that 
one has the standard estimate
\[
\|\Gamma_\varphi-\Gamma_\psi\|_1\leq h(\varphi-\psi)^{1/2}(h(\varphi)^{1/2}+h(\psi)^{1/2})
\;\;\;.
\]
Therefore one may extend $\Gamma$ to a quadratic form with domain $D(\overline h)$ and with values in $L_1(\Ri^d)$,
by continuity.
The extended form then defines an extension of the  bilinear form $\Gamma$  to $D(\overline h)\times D(\overline h)$
by polarization.
Since the coefficients $c_{ij}$ are  bounded the representation (\ref{econe2}) is valid if 
$\varphi\in  W^{1,2}(\Ri^d)\subseteq D(\overline h)$ but not in general.

If $\varphi,\psi\in C_c^\infty(\Ri^d)$ then it follows from 
(\ref{econe2}) by  Leibniz' rule that
\begin{equation}
\Gamma_{\varphi\psi}=\varphi^2\,\Gamma_{\psi}+2\,\varphi\,\psi \,\Gamma_{\varphi,\psi}+\psi^2\,\Gamma_{\varphi}
\;\;\;.
\label{econe2.1}
\end{equation}
Therefore
\begin{equation}
\Gamma_{\varphi\psi}\leq 2\,\varphi^2\,\Gamma_{\psi}+2\,\psi^2\,\Gamma_{\varphi}
\label{econe2.2}
\end{equation}
by the Cauchy-Schwarz inequality.
But it follows from the theory of Dirichlet forms that $D(\overline h)\cap L_\infty(\Ri^d)$
is an algebra.
Therefore one easily deduces that (\ref{econe2.2}) extends to all 
$\varphi,\psi\in D(\overline h)\cap L_\infty(\Ri^d)$.
Note that integration of (\ref{econe2.2}) gives an alternative version of the inequality in Lemma~\ref{lcdeg2.01}.
This is crucial in the following proofs.

Now we turn to the proof of Theorem~\ref{tcdeg1.1}.

\smallskip
\noindent{\bf Proof of Theorem~\ref{tcdeg1.1}}$\;$
\ref{tcdeg1.1-1}$\Rightarrow$\ref{tcdeg1.1-2}$\;$
First, fix a point $y\in A$ at which the surface is differentiable
and  choose $r$ sufficiently small that the ball $B=B(y\,;r)=\{x\in\Ri^d:|x-y|<r\}$ is bisected by
$A$ into two disjoint open subsets $B_+$ and $B_-$.

Secondly, fix $\Phi\in C_c^\infty(B)$ with $0\leq \Phi\leq 1$ and $\Phi(y)=1$.
Then define $h_\Phi$ by
\[
h_\Phi(\varphi)=\sum^d_{i,j=1}(\partial_i\varphi,\Phi \,c_{ij}\partial_j\varphi)
=\int_{\Ri^d} dx\,\Phi(x)\,\Gamma_\varphi(x)
\]
with domain $D(h_\Phi)=C_c^\infty(\Ri^d)$.
Then $h_\Phi$ corresponds to the divergence form operator with Lipschitz coefficients $\Phi\, c_{ij}$
 and is consequently closable.
Moreover $C^\infty(B)$ is a core of $h_\Phi$.
The closure $\overline h_\Phi$ of $h_\Phi$ is a Dirichlet form.
Let $H_\Phi$ and $S^\Phi$ denote the corresponding operator and semigroup.
Since $\Phi\in C_c^\infty(B)$  one can also view $h_\Phi$ as a form on $L_2(B)$
and the operator and semigroup as acting on the  subspace $L_2(B)$ of $L_2(\Ri^d)$,
i.e., $L_2(\Ri^d)$ has a canonical decomposition  $L_2(B)\oplus L_2(B^{\rm c})$
and the operator and semigroup act on the first component.
 Moreover, $\one_B\in D(h_\Phi)$ and $h_\Phi(\one_B)=0$.
Therefore $S^\Phi$ is conservative on $L_2(B)$, i.e., $S^\Phi_t\one_B=\one_B$ for all $t>0$.

Thirdly,  $0\leq C_h(A\cap B)\leq C_h(A)=0$ by positivity and monotonicity of the capacity and by Condition~\ref{tcdeg1.1-1}.
But if $U\subset B$ is open and $U\supset A\cap B$ then
\begin{eqnarray*}
C_h(U)&=&\inf\{\,\overline h(\varphi)+\|\varphi\|_2^2:\varphi\in D(\overline h)\,,\, 
\varphi\geq1\;{\rm on }\; U\,\}\\[8pt]
&\geq &\inf\{\,\overline h_\Phi(\varphi)+\|\varphi\|_2^2:\varphi\in D(\overline h_\Phi)\,,\, 
\varphi\geq1\;{\rm on }\; U\,\}=C_{h_\Phi}(U)
\;\;\;.
\end{eqnarray*}
Therefore $0\leq C_{h_\Phi}(A\cap B)\leq C_h(A\cap B)=0$.

Fourthly, since $C_{h_\Phi}(A\cap B)=0$ it follows from the proof of Proposition~\ref{pcdeg2.1} that 
$C^\infty(B \backslash(A\cap B))$ is a core of $h_\Phi$.
But each $\varphi\in C^\infty(B \backslash(A\cap B))$  has a unique decomposition
$\varphi=\varphi_++\varphi_-$ with $\varphi_\pm\in C^\infty(B_\pm)$.
Specifically $\varphi_+=\varphi$ on $B_+$ and $\varphi_+=0$ on $B_-$.
Note that $\varphi_\pm=\varphi\one_{B_\pm}$. 
Moreover,\begin{equation}
h_\Phi(\varphi)=h(\varphi_+)+h(\varphi_-)
=h_\Phi(\varphi\one_{B_+}) +h_\Phi(\varphi\one_{B_-})
\label{econe4}
\end{equation}
for all $\varphi$ in the core $C_c^\infty(B \backslash(A\cap B))$  of $h_\Phi$ because
the components $\varphi_+$ and $\varphi_-$ are $C^\infty$-functions
with disjoint supports.

Fifthly, if $\varphi\in D(\overline h_\Phi)$ one may choose $\varphi_n\in C^\infty(B \backslash(A\cap B))$
 such that $\overline h_\Phi(\varphi-\varphi_n)\to0$ and $\|\varphi-\varphi_n\|_2^2\to0$ as $n\to\infty$.
But 
\[
h_\Phi(\varphi_{n,+}-\varphi_{m,+})\leq h_\Phi(\varphi_n-\varphi_m)\;\;\;\;\;{\rm and }
\;\;\;\;\;\|\varphi_{n,+}-\varphi_{m,+}\|_2\leq \|\varphi_n-\varphi_m\|_2\;\;\;.
\]
Moreover, $\|\varphi_{n,+}-\varphi\one_{B_+}\|_2\to0$ as $n\to\infty$.
Therefore $\varphi\one_{B_+}\in D(\overline h_\Phi)$ and 
$\overline h_\Phi(\varphi\one_{B_+})=\lim_{n\to\infty}h_\Phi(\varphi_{n,+})$.
Similarly $\varphi\one_{B_-}\in D(\overline h_\Phi)$ and 
$\overline h_\Phi(\varphi\one_{B_-})=\lim_{n\to\infty}h_\Phi(\varphi_{n,-})$.
Therefore
\begin{equation}
\overline h_\Phi(\varphi)=\overline h_\Phi(\varphi\one_{B_+})+\overline h_\Phi(\varphi\one_{B_-})
\label{econe4.1}
\end{equation}
by taking limits of (\ref{econe4}) with $\varphi$ replaced by $\varphi_n$.

Sixthly, it follows from (\ref{econe4.1}) and Lemma~6.3 of \cite{ERSZ1} that $S^\Phi_tL_2(B_\pm)\subseteq L_2(B_\pm)$
for all $t>0$.
Therefore if $\varphi\in C_c^\infty(B)$ then
\[
(\one_{B_+},\varphi)=
(S^\Phi_t\one_{B_+},\varphi)=
(\one_{B_+},S^\Phi_t\varphi)
\]
for all $t\geq0$.
Hence
\[
0=(\one_{B_+},\varphi)-(\one_{B_+},S^\Phi_t\varphi)
=\int^t_0ds\,(\one_{B_+},S^\Phi_sH_\Phi\varphi)
=t\,(\one_{B_+},H_\Phi\varphi)
\]
for all $t>0$.
But then
\[
0=(\one_{B_+},H_\Phi\varphi)=
\sum^d_{i,j=1}\int_{B_+} dx\,(\partial_i\,\Phi\,c_{ij}\partial_j\,\varphi)(x)
=\int_{B_+} dx\,\divv\Psi
\]
where $\Psi=\sum^d_{j=1}\Phi\,c_{ij}\,\partial_j\varphi$.
Next one deduces from   Stokes' theorem that
\[
\sum^d_{i,j=1}\int_{\partial{B_+}} dS\,(n_{B_+,i}\,\Phi\,c_{ij})(x)(\partial_j\,\varphi)(x)=0
\;\;\;.
\]
Since this is valid for all $\varphi\in C_c^\infty(B)$ it follows that 
$\sum^d_{i=1}\Phi\,n_{A\cap B,i}\,c_{ij}=0$ on $A\cap B$.
In particular, since $\Phi(y)=1$ one has $\sum^d_{i=1}(n_{A\cap B,i}\,c_{ij})(y)=0$.
But the point $y$ was an arbitrary point of differentiability of the surface $A$ so the
flux must be zero at all such points.
Taking the scalar product with $n_{A\cap B}$ one obtains Condition~\ref{tcdeg1.1-2}.

\smallskip

\noindent\ref{tcdeg1.1-2}$\Rightarrow$\ref{tcdeg1.1-1}$\;$
The capacity satisfies $C_h(A\cup B)\leq C_h(A)+C_h(B)$ for all pairs of subsets $A$ and $B$.
Therefore it suffices to prove 
 $C_h(A_m)=0$ for all small subsets $A_m\subset A$,
i.e., the problem can again be localized.

Fix a point in $A$ at which the surface is differentiable and for
convenience let this point be the origin of coordinates.
Then let $B=B(0\,;r)$ be a small ball such that $A\cap B$ is differentiable.
If $r'\in \langle0,r\rangle$ we also set $B'=B(0\,;r')$.
Next choose local coordinates $(\tilde x_1,\ldots,\tilde x_d)$
in a neighbourhood of the origin 
such that the section $A\cap B$ is contained in the hypersurface
$\tilde x_1=0$.
The normal to  $A\cap B$ then 
corresponds to $(1,0,\ldots,0)$ in the new coordinates.

Let $\tilde B$ and $\tilde B'$ denote the images of $B$ and $B'$ respectively.
Then if $\varphi\in C_c^\infty(B)$ define $\widetilde\varphi\in C_c^\infty(\tilde B)$ by 
$\widetilde\varphi(\tilde x)=\varphi(x)$.
Note that 
\[
\|\varphi\|_2^2=\int_{\tilde B}d\tilde x\,J(\tilde x)\,|\widetilde\varphi(\tilde x)|^2
\]
and
\[
h(\varphi)=\sum^d_{i,j=1}\int_{\tilde B}d\tilde x\,J(\tilde x)\,\tilde c_{ij}(\tilde x)\,
(\tilde\partial_i\widetilde\varphi)(\tilde x)\,
(\tilde\partial_j\widetilde\varphi)(\tilde x)
\]
where $J$ denotes the Jacobian of the coordinate transformation, 
$\tilde\partial_i=\partial/\partial\tilde x_i$ and the coefficients $\widetilde C=(\tilde c_{ij})$ 
in the new coordinates are again symmetric, positive definite and Lipschitz continuous.
Moreover, Condition~\ref{tcdeg1.1-2} implies that $\tilde c_{11}(0,\tilde x_2,\ldots,\tilde x_d)=0$
and  the Lipschitz property implies $\tilde c_{11}(\tilde x_1,\tilde x_2\ldots,\tilde x_d)\leq a\,|\tilde x_1|$
for some $a>0$  within  the  neighbourhood $\tilde B$.

Next if $\varphi\in C_c^\infty(B)$ and $\varphi\geq 1$ on $A\cap B'$ we will construct approximants
$\varphi_n\in C_c^\infty(B)$ such that $\varphi_n\geq 1$ on $A\cap B'$ and $h(\varphi_n)+\|\varphi_n\|_2^2\to0$
as $n\to\infty$. 
This establishes that $C_h(A\cap B')=0$ and since this is valid for all $y\in A\cap B$ and all small balls
it follows that $C_h(A)=0$

The construction of the  approximants $\varphi_n$ to $\varphi$ follows the procedure used
in the proof of Proposition~6.5 of \cite{ERSZ1}.
First define  $\chi_n \colon \Ri \to [0,1]$ by
\begin{equation}
\chi_n(x)
= \left\{ \begin{array}{ll}
  0 & \mbox{if } x\leq -1 \;\;\; ,  \\[5pt]
  {-\log |x|}/{\log n} & \mbox{if } x\in\langle-1, -n^{-1}\rangle \;\;\; , \\[5pt]
  1 & \mbox{if } x \geq -n^{-1} \;\;\; .\label{ecdeg2.137}
         \end{array} \right.
\end{equation}
Note that $\chi_n$ is absolutely continuous and increasing and that
 $\lim_{n \to \infty} \chi_n = \one_{[0,\infty\rangle}$ pointwise.
In addition
\begin{eqnarray*}
\int_{\tilde B} d\tilde x\,J(\tilde x)\,\tilde c_{11}(\tilde x)\,(\chi_n'(\tilde x_1))^2
&\leq& a\,(\log n)^{-2}
\int^{-n^{-1}}_{-1}d\tilde x_1\,
\int_{\tilde B}d\tilde x_2\ldots d\tilde x_d\,
J(\tilde x)\,|\tilde x_1|\,
|\tilde x_1|^{-2}\\[5pt]
&\leq& a'\, (\log n)^{-1}
\end{eqnarray*}
for all $n \in \Ni$.
Therefore
\begin{equation}
\lim_{n \to \infty}\int_{\tilde B} d\tilde x\,J(\tilde x)\,\tilde c_{11}(\tilde x)\,(\chi_n'(\tilde
x_1))^2=0\;\;\;.\label{ecdeg2.99}
\end{equation}
Now introduce $\xi_n$  by $\xi_n(x)=\chi_n(x)\wedge \chi_n(-x)$.
Then one also has 
\begin{equation}
\lim_{n \to \infty}\int_{\tilde B} d\tilde x\,J(\tilde x)\,\tilde c_{11}(\tilde x)\,(\xi_n'(\tilde
x_1))^2=0\;\;\;.\label{ecdeg2.100}
\end{equation}
The approximants  $\varphi_n\in C_c^\infty(B)$ are now defined  such that 
$\widetilde\varphi_n(\tilde x)=\widetilde\varphi(\tilde x)\xi_n(\tilde x_1)$.

   It follows from this construction that 
\[
\|\varphi_n\|_2^2=\int_{\tilde B}d\tilde x\,J(\tilde x)\,
\xi_n(\tilde x_1)^2\,|\widetilde\varphi(\tilde x)|^2
=\int d\tilde x_1\,\xi_n(\tilde x_1)^2\,\psi(\tilde x_1)
\]
where $\psi$ is bounded with compact support.
Therefore $\|\varphi_n\|_2\to0$ as $n\to\infty$ because the $\xi_n$ converge
almost everywhere to zero.
Next, however,
\[
\Gamma_{\widetilde\varphi_n}=
\sum^d_{i,j=1}\tilde c_{ij}\,(\tilde\partial_i\widetilde \varphi_n)
\,(\tilde \partial_j\widetilde \varphi_n)
\leq 
2\,\tilde c_{11}\,(\xi'_n)^2\,\widetilde \varphi^2+2\,
\xi_n^2\,\Gamma_{\widetilde \varphi}
\]
where we have used (\ref{econe2.2}).
Therefore
\[
h(\varphi_n)=\int_{\tilde B}\, J\,\Gamma_{\widetilde \varphi_n}
\leq 2\,\Big(\int_{\tilde B} \,J\,\tilde c_{11}\,(\xi'_n)^2\Big)\, \|\widetilde\varphi\|^2_\infty
+2\int_{\tilde B} \,J\,\Gamma_{\widetilde \varphi}\,\xi_n^2
\;\;\;.
\]
But the first term on the right hand side tends to zero as $n\to\infty$ by (\ref{ecdeg2.100}) and the second
tends to zero because $J\,\Gamma_{\widetilde \varphi}$ is integrable and the $\xi_n$ converge almost everywhere to zero.
Therefore $C_h(A\cap B')=0$ and consequently $C_h(A)=0$.
\hfill$\Box$

\ruimte

The foregoing local estimates can be used to establish a strictly positive lower bound 
on the capacity of hypersurfaces for weakly subelliptic operators.
It is necessary to assume Lipschitz continuity of the surface but continuity of the 
coefficients is not necessary.

\begin{prop}\label{pcdeg3.1}
Assume the form $h$ with $L_\infty$-coefficients $c_{ij}$ satisfies the subellipticity estimate
\begin{equation}
h(\varphi)\geq \mu\,\|\Delta^{(1-\delta)/2}\varphi\|_2^2-\nu\,\|\varphi\|_2^2
\label{ecdeg3.69}
\end{equation}
 for some $\mu>0,\nu\geq 0$, $\delta\in[0,1\rangle$  and all $\varphi\in C_c^\infty(\Ri^d)$.
Let $h_0$ denote the corresponding viscosity form.
Further let $A$ be a Lipschitz continuous hypersurface in~$\Ri^d$.
Then 
\[
C_{h_0}(A)>0
\]
whenever $\delta\in[0,1/2\rangle$.
\end{prop}
\proof\
Again fix a point in $A$ at which the surface is differentiable and let this point be the origin of coordinates.
Then let $B=B(0\,;r)$ be a small ball such that $A\cap B$ is differentiable
and  choose local coordinates $(\tilde x_1,\ldots,\tilde x_d)$
in a neighbourhood of the origin 
such that the section $A\cap B$ is contained in the hypersurface
$\tilde x_1=0$.
It suffices to prove that $C_{h_0}(A\cap B)>0$.
Set $D=\{\,\tilde x: \tilde x_1=0\,, \,|\tilde x_2|^2+\ldots +|\tilde x_d|^2<\rho^2\,\}$.
If $\rho$ is sufficiently small $D\subset A\cap B$ and it suffices to prove $C_{h_0}(D)>0$.
But if $U_\varepsilon =\langle-\varepsilon, \varepsilon\rangle\times D$ it now suffices to prove
that $\lim_{\varepsilon\to0}C_{h_0}(U_\varepsilon)>0$.
One can, however, estimate $C_{h_0}(U_\varepsilon)>0$ with the aid of (\ref{ecdeg2.73}).

It is convenient to evaluate (\ref{ecdeg2.73}) in the new coordinates.
But since the local transformation of coordinates is non-singular and the Laplacian is unchanged
up to equivalence one has
\begin{eqnarray*}
(\one_{U_\varepsilon},(I+\Delta)^{-(1-\delta)}\one_{U_\varepsilon})&\leq &c\int_{\Ri^d}dx\,\one_{U_\varepsilon}(x)
((I+\Delta)^{-(1-\delta)}\one_{U_\varepsilon})(x)\\[5pt]
&=&c\int_{\Ri^d}dp\,(\tilde \one_{U_\varepsilon}(p))^2(1+p^2)^{-(1-\delta)}
\end{eqnarray*}
by Fourier transformation.
Since $\one_{U_\varepsilon}=\one_{\langle-\varepsilon, \varepsilon\rangle}\,\one_D$ it then follows that 
\begin{eqnarray*}
\hspace{-5pt}(\one_{U_\varepsilon},(I+\Delta)^{-(1-\delta)}\one_{U_\varepsilon})&\leq &
4\,\varepsilon^2\,c\int_{\Ri^{(d-1)}}d\tilde p\,(\tilde \one_D(\tilde p))^2
\int_\Ri dp_1\,(1+\tilde p^2+p_1^2)^{-(1-\delta)}(\sin(\varepsilon p_1)/(\varepsilon p_1))^2\\[5pt]
&\leq& 4\,\varepsilon^2 \,c_\delta\int_{\Ri^{(d-1)}}d\tilde p\,(\tilde \one_D(\tilde p))^2=4\,\varepsilon^2 \,c_\delta\,|D|
\end{eqnarray*}
where we have set $\tilde p=(p_2,\ldots,p_d)$.
Note that $c_\delta<\infty$ if, and only if, $\delta\in[0,1/2\rangle$.
Therefore it follows from (\ref{ecdeg2.73})   that
\[
C_{h_0}(U_\varepsilon)\geq 4\,\omega_\delta\,\varepsilon^2\,|D|^2\,(\one_{U_\varepsilon},(I+\Delta)^{-(1-\delta)}\one_{U_\varepsilon})^{-1}
\geq (\omega_\delta/c_\delta)\,|D|
\]
uniformly for $\varepsilon>0$.
\hfill$\Box$

\ruimte

 Proposition~\ref{pcdeg3.1} and Theorem~\ref{tcdeg1.2} have an immediate corollary.
If the coefficients $c_{ij}$ are Lipschitz continuous and if $h$ satisfies the subellipticity condition 
(\ref{ecdeg3.69}) with $\delta<1/2$ then it is impossible to have separation into independent subsystems with 
Lipschitz boundaries.

In Theorems \ref{tcdeg1.2} and \ref{tcdeg1.1} there are no restrictions on
the connectedness properties of the hypersurface of separation.
In fact one can construct examples in which it is disconnected or multiply
connected.
One interesting situation occurs for periodic coefficients.
We illustrate this with another one-dimensional example.
One can, however, construct examples in higher dimensions and strict periodicity
is not essential.
\begin{exam}\label{exldeg4.1}
Let $d=1$. Then $h(\varphi)=\int_\Ri c\,(\varphi')^2$ with $c\geq0$.
Choose $c$ such that $c(x)=(1-\cos2\pi x)^\delta$ with $\delta>0$.
Then $c$ has zeros of order $2\delta$ at the integer points.

First, if $\delta\geq 1/2$ then $C_{\overline h}(\{n\})=0$ for each $n\in\Ni$.
So separation occurs at each zero.
This follows because the function $\xi_n$ used in the proof of Theorem~\ref{tcdeg1.1}
 satisfies $\xi_n=1$ on $\langle-1/n,1/n\rangle$ and $h(\xi_n)+\|\xi_n\|_2^2\to0$
as $n\to\infty$.

Secondly if $\delta<1/2$ we argue that $C_{\overline h}(\{x\})>0$  for each $x\in \Ri$.
Hence there is no separation.
We begin by observing that
\[
h(\varphi)+\|\varphi\|_2^2\geq (E\varphi',\,c\,E\varphi')+(E\varphi,E\varphi)
\]
where $E$ is the orthogonal projection from $L_2(\Ri)$ onto $L_2(-2/3, 2/3)$
and $\varphi\in C_c^\infty(\Ri)$.
Now on the interval $\langle-2/3,2/3\rangle$ one has $c(x)\geq a\,(x^2/(1+x^2))^\delta$
for some $a>0$.
Since $\delta<1/2$ it follows, however, from Theorem~3.6 in \cite{Stri} that 
$\Delta^\delta\geq \sigma |x|^{-2\delta}$ where  $\Delta$ is again the self-adjoint 
version of $-d^2/dx^2$ on $L_2(\Ri)$, i.e., the Laplacian.
One then has a bound $a\,(x^2/(1+x^2))^\delta\geq a_0\,(I+\Delta)^{-\delta}$ on $L_2(\Ri)$.
Next let $\Delta_D$  denote the Laplacian on $L_2(\Ri)$ with Dirichlet  boundary conditions
 at the points $x=\pm 2/3$
Therefore
\[
(E\varphi',\,c\,E\varphi')\geq a_0\,(E\varphi',(I+\Delta)^{-\delta}E\varphi')
\geq a_0\,(E\varphi',(I+\Delta_D)^{-\delta}E\varphi')
\;\;\;.
\]
But $E\,(I+\Delta_D)^{-\delta}E=(I+\tilde\Delta_D)^{-\delta}$ where 
$\tilde\Delta_D$ is the restriction of $\Delta_D$ to 
$L_2(-2/3, 2/3)$, i.e., the operator with Dirichlet boundary conditions at the endpoints.
But 
\begin{eqnarray*}
(\varphi',(I+\tilde \Delta_D)^{-\delta}\varphi')_I
&=&(d\varphi,(I+\tilde \Delta_D)^{-\delta}d\varphi)_I\\[5pt]
&=&(d\varphi,d(I+\tilde \Delta_N)^{-\delta}\varphi)_I
=(\tilde\Delta_N^{1/2}(I+\tilde \Delta_N)^{-\delta/2}\varphi,
\tilde\Delta_N^{1/2}(I+\tilde \Delta_N)^{-\delta/2}\varphi)_I
\end{eqnarray*}
where the scalar product is on $L_2(-2/3, 2/3)$, $d$ denotes the closed operator 
of differentiation  with no boundary conditions and  $\tilde\Delta_N$ 
is the Laplacian with Neumann boundary conditions at the endpoints.
It then follows easily from combination of these estimates that one has a subellipticity bound
\[
h(\varphi)+\|\varphi\|_2^2\geq a_1
((I+\tilde\Delta_N)^{(1-\delta)/2}\varphi,
(I+\tilde\Delta_N)^{(1-\delta)/2}\varphi)_I
\]
with $a_1>0$.
Therefore if $x\in[-1/2,1/2]$ and $U\subset \langle-2/3,2/3\rangle$ is an open subinterval
containing $x$ with length $|U|$ then one estimates as in the proof of Proposition~\ref{pcdeg3.1}
that
\[
C_{\overline h}(U)\geq a_1\,|U|^2(\one_U,(I+\tilde\Delta_N)^{-(1-\delta)/2}\one_U)_I^{-1}
\;\;\;.
\]
But $\tilde\Delta_N$ has zero as smallest eigenvalue with constant eigenfunction.
Therefore one has $C_{\overline h}(U)\geq 3a_1/4>0$.
Hence $C_{\overline h}(\{x\})>0$.
The result then extends to all $x\in\Ri$ by periodicity.

In this example $S\cup L_\infty$ is irreducible on $L_2(\Ri)$ if $\delta<1/2$ but it is not irreducible
if $\delta\geq 1/2$.
In the latter case the space $L_2(\Ri)$ can be decomposed as a direct sum of copies of $L_2(0,1)$ each
of which is invariant under $S\cup L_\infty$. 
The family of operators acts irreducibly on each component space.
\end{exam}

\section{Concluding remarks}\label{Scdeg3}

In this section we briefly discuss various aspects of 
Theorems~\ref{tcdeg1.1} and \ref{tcdeg1.2} and possible extensions.

Although both theorems were based on the assumption 
of Lipschitz continuity of  the coefficients $c_{ij}$ the only essential use
of this  property was in an open neighbourhood  $U$ of the hypersurface $A$.
The continuity was used in an inessential manner to deduce closability of the form $h$ 
through  the theory of the Friederichs extension. 
But this can be avoided. 
The form is also closable if one has Lipschitz continuity of the coefficients on $U$ 
and strong ellipticity on $\Ri^d\backslash U$.
This follows straightforwardly by decomposing $h$ into two positive components
 on $U$ and $\Ri^d\backslash U$ respectively.
Each component is then closable and this suffices to deduce closability of the sum.

The Lipschitz continuity near the hypersurface is not necessary for separation.
This is  illustrated in one-dimension by Example~\ref{excdeg2.1}.
Then there is only one coefficient $c$ with  $c(x)=O(|x|)$ as $x\to0+$ and 
 $c(x)=O(|x|^\delta)$ with $\delta\in[0,1\rangle$ as $x\to0-$.
But 
one has separation into two components, the left and right half-lines, by the arguments
of Proposition~6.5 of \cite{ERSZ1}.
In higher dimensions one can also have a similar effect of separation corresponding
to degeneracy on one side of the hypersurface $A$ and this can be complicated
as the degeneracy may vary from side to side.

Next we comment on several properties of the energy density $\Gamma$.
First, we note that the piecewise differentiable surface $A$ is characterized 
as the zero set of a Lipschitz function $\eta$.
Therefore one can extend the energy density to $\eta$ by setting
\[
\Gamma_\eta=\sum^d_{i,j=1}c_{ij}(\partial_i\eta)(\partial_j\eta)
\;\;\;.
\]
Then the zero flux property, $(n_A,Cn_A)=0$ on $A$, corresponds to the property $\Gamma_\eta=0$
on~$A$.
This in turn shows that the surface of separation is a subset of the surface $\mu_C=0$ where
$\mu_C$ is the lowest eigenvalue of the coefficient matrix $C$.

Secondly, the energy density $\Gamma$ is directly related to various definitions of distance.
For example the Riemannian distance corresponding to $H$ can be defined by a shortest path
algorithm or by setting
\begin{equation}
d(x\,;y)=\sup_{\psi\in D}|\psi(x)-\psi(y)|
\label{ecdeg3.1}
\end{equation}
where the set of variational functions is defined by
\begin{equation}
D=\{\,\psi\in C_c^\infty(\Ri^d): \|\Gamma_\psi\|_\infty\leq 1\,\}
\;\;\;.
\label{ecdeg3.2}
\end{equation}
This distance is automatically continuous and is not ideally suited to the discussion
of separation phenomena and discontinuous behaviour.
The natural solution is to modify this definition to allow a larger class $D_0$
of  functions including some discontinuous ones.
But the introduction of discontinuous functions  requires  a modification of (\ref{ecdeg3.1}).
One approach is to restrict consideration to  the distance between  measurable subsets.
 We refer to \cite{ERSZ2} for details 
(see also \cite{Stu2} \cite{HiR}).
The natural class of variational functions in this extended definition is given by
$D_0=\{\psi\in D^\infty(h):\|\Gamma_\psi\|_\infty\leq1\}$
where the set $D^\infty(h)$ is given by $D^\infty(h)=\{\psi\in D(\overline h)\cap L_\infty:\|\Gamma_\psi\|_\infty<\infty\}$.
In fact the subset $D_0$ is not adequate for the discussion of unbounded subsets  and 
global properties because it consists of functions which are zero at infinity;
an even larger set is required for these purposes (see \cite{ERSZ2}).
Nevertheless $D_0$ should be adequate for local properties (see  \cite{HiR}).
Note that if one adopts these definitions then it is evident that if the separation phenomena
analyzed in Section~\ref{Scdeg2} occur then $D^\infty(h)\not\subseteq C_0(\Ri^d)$.
Indeed in the proof of Theorem~\ref{tcdeg1.1} we established that for each 
$\varphi\in C_c(\Ri^d)$ one has 
$\varphi\one_{B_+}\in D^\infty(h)$.
But if the support of $\varphi$ intersects the interface between $B_+$ and $B_-$ then
$\varphi\one_{B_+}$ is discontinuous.
It is also likely that a local converse  is true, i.e.,
it is possible that   for Lipschitz continuous coefficients local separation
occurs if and only if  $D^\infty(h)\not\subseteq C_0(\Ri^d)$.
This is the situation for the form $h_\delta$ in Example~\ref{excdeg2.1}.
If $\delta<1/2$ then $D(\overline h_\delta)\subset C_0(\Ri)$.
Hence $D^\infty(h_\delta)\subseteq C_0(\Ri)$.
But if $\delta\in[1/2,1\rangle$ separation occurs and 
 $D^\infty(h_\delta)\not\subseteq C_0(\Ri)$

Finally we note that if $D^\infty(h)\subseteq C_0(\Ri^d)$ then each open
subset $\Omega$  has boundary $\partial \Omega$ with 
$C_h(\partial\Omega)>0$.
Indeed if $C_h(\partial\Omega)=0$ then one has separation by Theorem~\ref{tcdeg1.0}
and this implies that 
 $D^\infty(h)\not\subseteq C_0(\Ri^d)$.

\section*{Acknowledgements}

This work was supported by the Australian Research Council's Discovery
Grant Program.
The authors are grateful to Tom ter Elst for critically reading the manuscript.

\end{document}